\documentclass[12pt,oneside,english]{amsart}
\usepackage[latin1]{inputenc}
\usepackage{amssymb}

\makeatletter

\newtheorem{lemma}{Lemma}
\newtheorem{example}{Example}
\newtheorem{definition}{Definition}

\numberwithin{equation}{section}
\usepackage{babel}

\makeatother
\begin{document}
\baselineskip=17pt

\title[Combinatorial minors ]{Combinatorial minors for matrix functions and their
applications}

\author{Vladimir Shevelev}
\address{Department of Mathematics \\Ben-Gurion University of the
 Negev\\Beer-Sheva 84105, Israel. e-mail:shevelev@bgu.ac.il}

\subjclass{15A15; key words and phrases: square matrices, permanent, cycle permanent, $\omega$-permanent, cycle index of permutations with restricted positions.}

\begin{abstract}
As well known, permanent of a square (0,1)-matrix $A$ of order $n$ enumerates the
permutations $\beta$ of $1,2,...,n$ with the incidence matrices $B\leq A.$
To obtain enumerative information on even and odd permutations with
condition $B\leq A,$  we should calculate two-fold vector $(a_1,a_2)$ with $a_1+a_2
=per A.$ More general, the introduced $\omega$-permanent, where
$\omega=e^{2\pi i/m},$ we calculate as $m$-fold vector. For these and other matrix
functions we generalize the Laplace theorem of their expansion over elements of the
first row, using the defined so-called "combinatorial minors". In particular,
 in this way, we calculate the cycle index of permutations with condition $B\leq A.$
\end{abstract}

\maketitle

\section{Introduction}

Let $\mathcal{S}^{(n)}$ be symmetric group of permutations of numbers $\{1,...,n\}.$ Let $A=\{a_{ij}\}$ be square matrix of order $n.$ Recall that the permanent of $A$ is defined by formula (\cite{4})
\begin{equation}\label{1.1}
per A=\sum_{s\in \mathcal{S}^{(n)}}\prod_{i=1}^n a_{i,s(i)}.
\end{equation}
If $A$ is a $(0,1)$ matrix, then it defines a class $\mathcal{B}=\mathcal{B}(A)$ of permutations with restricted positions, such that the positions of its zeros are prohibited. Such class could be equivalently defined by a simple inequality: a permutation $\pi\in \mathcal{B}$ if and only if for its incidence matrix $P$ we have $P\leq A.$ One of the most important application of $per A$ consists of the equality $|\mathcal{B}|=per A.$ Thus $per A$ enumerates
permutations with restricted positions of the class  $\mathcal{B}(A).$ \newline
\indent Let $\gamma(\pi)$ be number of independent cycles of $\pi,$ including cycles of length 1. Then the difference
$d(\pi)=n-\gamma(\pi)$ is called \slshape decrement \upshape \enskip of $\pi$ (\cite{3}).
 Permutation $\pi$ is called even (odd), \enskip if $d(\pi)$ is even (odd). Note that the determinant of matrix $A$ could be defined by formula
 \begin{equation}\label{1.2}
det A=\sum_{even\enskip s\in \mathcal{S}^{(n)}}\prod_{i=1}^n a_{i,s(i)}-\sum_{odd\enskip s\in \mathcal{S}^{(n)}}\prod_{i=1}^na_{i,s(i)}.
\end{equation}
Since, evidently, we also have
 \begin{equation}\label{1.3}
per A=\sum_{even\enskip s\in \mathcal{S}^{(n)}}\prod_{i=1}^n a_{i,s(i)}+\sum_{odd\enskip s\in \mathcal{S}^{(n)}}\prod_{i=1}^n a_{i,s(i)},
\end{equation}
then the numbers of even and odd permutations of class $\mathcal{B}(A)$ are given by vector
\begin{equation}\label{1.4}
(\frac{1}{2}(per A+det A),\enskip \frac{1}{2}(per A-det A)).
\end{equation}
Note that, in contrast to permanent, there exist methods of very fast calculation of $det A.$ Therefore, the enumerative information given by (\ref{1.4}) one can obtain approximately for the same time as the number $|\mathcal{B}(A)|$ given by (\ref{1.1}).\newline
\indent Let $m\geq3$ and $0\leq k<m$ be given integers. We say that a permutation $\pi$ belongs to \slshape class $k$ modulo \upshape $m$ $(\pi\in \mathcal{S}_{k, m}^{(n)}),$ if $d(\pi)\equiv k\pmod m.$  The first problem under our consideration in this paper is the enumeration of permutations $\pi\in \mathcal{B(A)}$ of class $k$ modulo $m.$ It is clear that this problem is a natural generalization of problem of enumeration of even and odd permutation with restricted positions which is solved by ({1.4}). In order to solve this more general problem, put $\omega=e^{\frac{2\pi i}{m}}$ and introduce a new matrix function which we call \slshape $\omega$-permanent.\upshape
\begin{definition}\label{d1} Let $A$ be a square matrix of order $n.$ We call $\omega$-permanent of matrix $A$ the following matrix function

$$per_\omega A=\sum_{\enskip s\in \mathcal{S}_{0,m}^{(n)}}\prod_{i=1}^n a_{i,s(i)}+\omega\sum_{\enskip s\in \mathcal{S}_{1,m}^{(n)}}\prod_{i=1}^n a_{i,s(i)}+$$
\begin{equation}\label{1.5}
\omega^2\sum_{\enskip s\in \mathcal{S}_{2,m}^{(n)}}\prod_{i=1}^n a_{i,s(i)}+...+\omega^{m-1}\sum_{\enskip s\in \mathcal{S}_{m-1,m}^{(n)}}\prod_{i=1}^n a_{i,s(i)}.
\end{equation}
\end{definition}
Note that, if $m=1,$ then $per_\omega A=per A,$ and if $m=2,$ then $per_\omega A=det A.$ In case $m\geq3,$ every sum in ({1.5}) essentially differs from permanent. Therefore, the known methods of evaluation of permanent (\cite{4}, ch. 7) are not applicable. However, using so-called "combinatorial minors", below we find an expansion $per_\omega A$ over the first row of matrix $A.$ This allows to reduce a problem of order $n$ to a few problems of order $n-1.$ \newline
\indent The second problem under our consideration is another important problem of enumeration of full cycles of length $n$ with restricted positions. In connection with this problem, we introduce another new matrix function which we call \slshape cyclic permanent. \upshape
\begin{definition}\label{d2} Let $A$ be square matrix of order $n.$ The number
\begin{equation}\label{1.6}
Cycl(A)=\sum_{s}\prod_{i=1}^na_{i,s(i)},
\end{equation}
where the summing is over all full cycles from $\mathcal{S}^{(n)},$ we call a cycle permanent of $A.$
\end{definition}
The third problem is a problems of enumeration of permutations with restricted positions with a restrictions on their cycle structure.\newline
\indent Denote $\gamma(s)$ the number of all cycles of permutation $s$ (including
cycles of length 1, if they exist). Recall (\cite {1}) that the absolute value of Stirling number $S(n,k)$ of the first kind equals to number of permutations $s\in\mathcal{S}^{(n)}$ with $\gamma(s)=k$ (\cite{1},\cite {5}). A natural generalization of Stirling numbers of the first kind is the following matrix function.
\begin{definition}\label{d3} The matrix function
\begin{equation}\label{1.7}
S(A; n,k)=\sum_{s\in\mathcal{S}^{(n)},\gamma(s)=k}\prod_{i=1}^na_{i,s(i)},
\end{equation}
where $n$ is order of square matrix $A,$ we call Stirling function of index $k.$
\end{definition}
Finally, recall (\cite{5}) that a permutation $s\in\mathcal{S}^{(n)}$ with $k_1$ cycles of length 1, $k_2$ cycles of length 2, and so on, is said to be of cycle structure $\overline{k}=(k_1, k_2,...,k_n).$ Denote $\nu(k_1,k_2,...,k_n)$ the number of permutations of class $\overline{k}=(k_1, k_2,...,k_n).$ Then the polynomial
\begin{equation}\label{1.8}
C(t_1, t_2,...,t_n)=\sum \nu(k_1,k_2,...,k_n)t_1^{k_1}t_2^{k_2}...t_n^{k_n}
\end{equation}
is called the cycle index of permutations of $\mathcal{S}^{(n)}.$ A natural generalization of the cycle index on permutations with restricted positions is
\begin{equation}\label{1.9}
C(A; t_1, t_2,...,t_n)=\sum \nu(A; k_1,k_2,...,k_n)t_1^{k_1}t_2^{k_2}...t_n^{k_n}
\end{equation}
where $\nu(A; k_1,k_2,...,k_n)$ is the number of permutations of class $\mathcal{B}(A)$ with the cycle structure $\overline{k}.$ In particular, if $A=J_n,$ where $J_n$ is $n\times n$ matrix
from 1's only, we have $C(A; t_1, t_2,...,t_n)=C( t_1, t_2,...,t_n). $
\section{Observations in case $m=2$ of $\omega$-permanent (determinant)}
 The case of determinant ($m=2$) of $\omega$-permanent is a unique case when it is easy to obtain a required enumerative information formally using formulas (\ref{1.4}). For the passage to a general case, it is important for us to understand how one can obtain such information from the definition (\ref{1.2}) of determinant only. Essentially, the required information is contained
 in vector
 \begin{equation}\label{2.1}
\overline{det}A=(\sum_{even\enskip s\in \mathcal{S}}\prod_{i=1}^n a_{i,s(i)}\enskip, \enskip -\sum_{odd\enskip s\in \mathcal{S}}\prod_{i=1}^n a_{i,s(i)}),
\end{equation}
\newpage
and it immediately disappears, if to add its components, or formally to use an identity of the form $1-1=0.$ Since $a_{i,j}=0$ or 1, then in sums of (\ref{2.1}) we do not use such
an identity. However, none of algorithms of fast calculation of determinant exists without
using of it. On the other hand, using the Laplace algorithm of expansion of
 determinant over
(the first) row, it can only be used at the last step. Therefore, if \slshape not to do \upshape \enskip the last step, we can obtain the required enumerative information.
\begin{example}\label{e1}
By the Laplace expansion, we have
$$ {det}\begin{pmatrix}0&1&1\\1&1&1\\1&1&1
\end{pmatrix}=-(1-1)+(1-1)=-1+1+1-1=2-2$$

and, if not to do the useless (with the enumerative point of view) last step, then we have
$$\overline{det}\begin{pmatrix}0&1&1\\1&1&1\\1&1&1
\end{pmatrix}=(2,-2).  $$
\end{example}
This means that there are two even and two odd permutations with the prohibited position $(1,1).$ \newline
 Of course, in general, the Laplace expansion of determinant over the first
 row in its classic form
 \begin{equation}\label{2.2}
det A=\sum_{j=1}^n (-1)^{j-1}a_{1,j}M_{1j},
\end{equation}
where $ M_{1j}$ is minor of element $a_{1j},$ i.e., determinant of the complementary to $a_{1j}$ submatrix $A_{1j},$ does not work for $\omega$-permanent.\newline
\indent Therefore, let us introduce for our aims a more suitable notion of so-called
 "combinatorial minor" of element $a_{ij}.$ Let the complementary to $a_{ij}$ submatrix $A_{ij}$ have the following $n-1$ columns
\begin{equation}\label{2.3}
c_1, c_2,..., c_{j-1}, c_{j+1},..., c_n.
\end{equation}
The first $j-1$ of these columns we change in the following cyclic order: $c_2, c_3,...,c_{j-1}, c_1.$ Then we obtain a new matrix ${A'}_{ij}$ with the columns
\begin{equation}\label{2.4}
 c_{2}, c_3,...,c_{j-1}, c_1, c_{j+1},..., c_n.
\end{equation}
Determinant of matrix ${A'}_{ij}$ we call \enskip \slshape combinatorial minor\upshape \enskip $(CM)_{ij}$ of element $a_{ij}$ . It is easy to see that
\begin{equation}\label{2.5}
(CM)_{i1} =M_{i1};\enskip (CM)_{ij} =(-1)^{j-2}M_{ij}, \enskip j=2,...,n.
\end{equation}
Therefore, e.g., expansion (\ref{2.2}) one can rewrite in the form
\begin{equation}\label{2.6}
det A=a_{1,1}(CM)_{11}-\sum_{j=2}^n a_{1,j}(CM)_{1j}.
\end{equation}
\newpage
In general, let us give a definition of combinatorial minors for arbitrary matrix function $X (A).$
\begin{definition}\label{d4} Let $X$ be matrix function defined on all square matrices of order $n\geq3.$ Let $A=\{a_{ij}\}$ be a square matrix of order $n$ and $A_{ij}$ be the complementary to $a_{ij}$ submatrix with columns $(\ref{2.3}).$ Denote ${A'}_{ij}$ a new square matrix of order $n-1$ with columns $(\ref{2.4}).$ Then the number
$X({A'}_{ij})$ is called a combinatorial minor of  $a_{ij}.$
\end{definition}

It appears that our observation (\ref{2.6}) has a general character. So, in Sections 4 we give a generalization of the Laplace expansion of type (\ref{2.6}) for $per_\omega A,$  $Cycl(A),$  $S(A; n, k)$ and $C_n(A; t_1, t_2,...,t_n).$

\section{Main lemma}

\begin{lemma}\label{L1}
Let $\pi\in \mathcal{S}^{(n)}$ with $\pi(j)=1$ and
 $\sigma=\sigma_j(\pi),\enskip j\geq2,$ such that

$$ \sigma(1)=\pi(2),\enskip \sigma(2)=\pi(3),\enskip ...,\enskip\sigma(j-2)=\pi(j-1),\enskip \sigma(j-1)=\pi(1),$$
 \begin{equation}\label{3.1}
   \sigma(j)=\pi(j)=1,\enskip\sigma(j+1)=\pi(j+1),\enskip...,\enskip\sigma(n)=\pi(n).
 \end{equation}
  Let, further, $\pi^*\in \mathcal{S}^{(n-1)}$ defined by the formula
\begin{equation}\label{3.2}
\pi^*(i)=\begin{cases}\sigma(i)-1,  \;\; if\;\;1\leq i\leq j-1,\\
\sigma(i+1)-1, \;\; if\;\;j\leq i\leq n-1.\end{cases}
\end{equation}
Then permutations $\pi$ and $\pi^*$ have the same number of cycles:
\begin{equation}\label{3.3}
\gamma(\pi)=\gamma(\pi^*).
\end{equation}
\end{lemma}
\bfseries Proof.\mdseries \enskip From (\ref{3.1})-(\ref{3.2}) we find
\begin{equation}\label{3.4}
\pi^*(i)=\begin{cases}\pi(i+1)-1,  \;\; if\;\;i\neq j-1,\\
\pi(1)-1, \;\; if\;\;i=j-1.\end{cases}
\end{equation}

Consider a cycle of $\pi$ containing element $j.$ Let it has length $l\geq2,$ such that
$$\pi(j)=1, \enskip\pi(1)=k_1,\enskip\pi(k_1)=k_2,\enskip...,\enskip\pi(k_{l-3})=k_{l-2},\enskip \pi(k_{l-2})=j. $$
 Beginning with the equality $\pi(1)=k_1,$ by (\ref{3.4}), this means that
$$ \pi^*(j-1)=k_1-1,\enskip \pi^*(k_1-1)=k_2-1,\enskip...,$$ $$ \pi^*(k_{l-3}-1)=k_{l-2}-1, \enskip\pi^*(k_{l-2}-1)=j-1. $$
Thus to a cycle of $\pi$ of length $l\geq2$ containing element $j$ corresponds a
cycle of length $l-1$ of $ \pi^*.$
Quite analogously, we verify that to cycle of $\pi$ of length $l\geq2$
not containing element $j$  corresponds a cycle of the same length of $ \pi^*.$
For example, to cycle of length $l\geq2$ of the form
\newpage
$$\pi(j-1)=k_1, \enskip\pi(k_1)=k_2,\enskip\pi(k_2)=k_3,\enskip...,\enskip\pi(k_{l-2})=k_{l-1},\enskip \pi(k_{l-1})=j-1 $$
(beginning with $\pi(k_1)=k_2),$ corresponds the cycle of the same length
$$ \pi^*(k_1-1)=k_2-1,\enskip \pi^*(k_2-1)=k_3-1,\enskip...,$$ $$ \pi^*(k_{l-2}-1)=k_{l-1}-1, \enskip\pi^*(k_{l-1}-1)=j-2 $$
such that $\pi^*(j-2)=\pi(j-1)-1=k_1-1.$\newline
$\blacksquare$\newline
\indent Note that the structure of Lemma \ref{L1} completely corresponds to the procedure of creating the combinatorial minors.

\section{Laplace expansions of type ({2.6}) of $per_\omega A,$ $Cycl(A),$ $S(A; n,k)$ and $C(A;t_1, t_2,...,t_n)$}
1) $per_\omega A.$  Consider all permutations $\pi$ with the condition $\pi(j)=1.$ Let $j$ corresponds to $j$-th column of matrix $A.$ Then the considered permutations correspond to diagonals of matrix $A$ having the common position $(1,j).$ If $j=1,$ then, removing the first row and column, we diminish on 1 the number of cycles of every such permutation, but also we diminish on 1 the number of elements of permutations. Therefore, the decrement of permutations does not change. If $j\geq2,$ consider continuation of these diagonals to the matrix of combinatorial minor ${A'}_{1,j}.$ Then, by Lemma \ref{L1}, the number of cycles of every its diagonal does not change and, consequently, the decrement is diminished by 1. This means that we have the following expansion of  $per_\omega A$ over the first row
\begin{equation}\label{4.1}
det_\omega A=a_{1,1}(CM)_{11}+\omega\sum_{j=2}^n a_{1,j}(CM)_{1j},
\end{equation}
where $(CM)_{1j},\enskip j\geq1,$ are combinatorial minors for $per_\omega A.$\newline
\indent Note that, as for determinant (see Section 2), in order to receive the required enumerative information, we should prohibit to use the identities of type $1+\omega+...+\omega^{m-1}=0.$\newline\newline
\indent 2)\enskip $Cycl(A).$  For $n>1,$ here we should ignore element $a_{11}.$ Consider all full cycles $\pi$ with the condition $\pi(j)=1.$ Let $j$ corresponds to $j$-th column of matrix $A.$ Then the considered full cycles correspond to diagonals of matrix $A$ having the common position $(1,j), \enskip j\geq2.$ Consider continuation of these diagonals to the matrix of combinatorial minor ${A'}_{1,j}.$ Then, by Lemma \ref{L1}, the number of cycles of every its diagonal does not change, i.e., they are full cycles of of order $n-1.$ Therefore, we have the following expansion of
\newpage
 $Cycle(A)$ over the first row of $A$
\begin{equation}\label{4.2}
Cycle(A)=\sum_{j=2}^n a_{1,j}(CM)_{1j},
\end{equation}
where $(CM)_{1j},\enskip j\geq1,$ are combinatorial minors for $Cycle(A).$\newline\newline
\indent 3) $S(A; n,k).$ From very close to 1) arguments, we have the following expansion of  $S(A; n,k)$ over the first row of $A$
\begin{equation}\label{4.3}
S(A; n,k)=a_{1,1}(CM)_{11}^{(k-1)}+\sum_{j=2}^n a_{1,j}(CM)_{1j}^{(k)},
\end{equation}
where $(CM)_{1j}^{(k)},\enskip j\geq1,$ are combinatorial minors of $S(A; n,k).$\newline
\indent Note that a close to ({4.3}) formula was found by the author in \cite{7} but using much more complicated way.\newline\newline
\indent 4) $C(t_1, t_2,...,t_n).$ We need lemma.
\begin{lemma}\label{L2}
Let
\begin{equation}\label{4.4}
\{a_{1j}, a_{k_{1}1}, a_{k_{2}k_{1}}, ... , a_{k_{r},k_{r-1}}, a_{j,k_{r}}\}
\end{equation}
be a cycle. Then
\begin{equation}\label{4.5}
\{a_{k_{1}1}, a_{k_{2}k_{1}}, ... , a_{k_{r},k_{r-1}}, a_{j,k_{r}}\}
\end{equation}
is a cycle with respect to the main diagonal of the matrix of combinatorial minor ${A'}_{1,j}.$
\end{lemma}
\bfseries Proof.\mdseries \enskip According to the construction of ${A'}_{1,j},$ the main its diagonal is
\begin{equation}\label{4.6}
\{a_{22}, a_{33},...,a_{j-1 j-1}, a_{j1}, a_{j+1 j+1},...,a_{nn}\}.
\end{equation}
With respect to this diagonal we have the following contour which shows that (\ref{4.5}) is, indeed, a cycle.
\begin{equation}\label{4.7}
\{a_{k_11}\rightarrow a_{j1}\rightarrow a_{jk_r}\rightarrow a_{k_rk_{r-1}}\rightarrow a_{k{r-1}k_{r-2}}\rightarrow...\rightarrow a_{k_2k_1} (\rightarrow a_{k_11})\}.
\end{equation}
$\blacksquare$\newline
\indent Quite analogously we can prove that to every another cycle of a diagonal containing element $a_{1j}$ correspond
the same cycle with respect to the main diagonal of the matrix ${A'}_{1,j}.$\newline
\indent Let $A$ be $(0,1)$ square matrix of order $n.$ Denote by $C^{(r)}(A; t_1, t_2,...,t_n)$ a partial cycle index of index (\ref{1.9}) of permutations $\pi\in \mathcal{B}(A)$ for which $\{1, \pi(1), \pi(2),...,\pi(r-1)\} $ is a cycle of length $r.$ Then we have
\begin{equation}\label{4.8}
\sum_{r=1}^n C^{(r)}(A; t_1, t_2,...,t_n)=C(A; t_1, t_2,...,t_n).
\end{equation}
\newpage
Therefore, it is sufficient to give an expansion of $C^{(r)}(A; t_1, t_2,...,t_n),\enskip r=1,...,n.$ First of all,
note that
\begin{equation}\label{4.9}
 C^{(1)}(A; t_1, t_2,...,t_n)=a_{11}t_1C(A; t_1, t_2,...,t_{n-1}).
\end{equation}
 Furthermore, using Lemmas \ref{L1}-\ref{L2}, we have
\begin{equation}\label{4.10}
C^{(r)}(A; t_1, t_2,...,t_n)=\frac{t_r} {t_{r-1}}\sum_{j=2}^{n}a_{1,j}(CM)_{1,j},
\end{equation}
where
\begin{equation}\label{4.11}
(CM)_{1,j}=C^{(r-1)}({A'}_{1,j}; t_1, t_2,...,t_{n-1}),\enskip j=2,...,n,
\end{equation}
are the combinatorial minors for cycle index of permutations with restricted positions.
Note that factor $\frac{t_r} {t_{r-1}}$ in (\ref{4.10}) corresponds to the diminution of the length of cycle (\ref{4.5}) with respect to length of  cycle (\ref{4.4}).\newline
\indent Thus formulas (\ref{4.8})-({4.11}) reduce the calculation of cycle index of $n$-permutations to the calculation of cycle index of $n-1$-permutations with a
rather simple computer realization of this procedure. Note that similar but much
more complicated procedure was indicated by the author in \cite{8}.

\section{An example of enumerating the permutations of classes 0,1,2 modulo 3 with
 restricted positions}
Let
$$ A=\begin{pmatrix}1&1&1&1&0\\0&1&0&1&1\\1&0&1&1&1\\1&1&1&0&0\\1&1&1&1&1
\end{pmatrix}. $$
Consider class $\mathcal{B}(A)$ of permutations with restricted positions and find the distribution of them over classes 0,1,2 modulo 3. We use $\omega$-permanent with $\omega=e^{\frac{2\pi i}{3}}$ and its expansion over elements of the first row, given by (\ref{4.1}). Recall that, for the receiving the required enumerative information, we should not use identities of type $1+\omega+\omega^2=0.$\newline
\indent We have
$$per_\omega A=per_{\omega} \begin{pmatrix}1&0&1&1\\0&1&1&1\\1&1&0&0\\1&1&1&1\end{pmatrix}+$$ $$\omega ( \enskip per_{\omega} \begin{pmatrix}0&0&1&1\\1&1&1&1\\1&1&0&0\\1&1&1&1\end{pmatrix}+per_{\omega} \begin{pmatrix}1&0&1&1\\0&1&1&1\\1&1&0&0\\1&1&1&1\end{pmatrix}+per_{\omega} \begin{pmatrix}1&0&0&1\\0&1&1&1\\1&1&1&0\\1&1&1&1\end{pmatrix})= $$
\newpage
$$per_{\omega} \begin{pmatrix}1&1&1\\1&0&0\\1&1&1\end{pmatrix}+\omega(\enskip per_{\omega}\begin{pmatrix}1&0&1\\1&1&0\\1&1&1\end{pmatrix}+ per_{\omega}\begin{pmatrix}1&1&0\\1&0&1\\1&1&1\end{pmatrix})+ $$
$$\omega^2(\enskip per_{\omega}\begin{pmatrix}1&1&1\\1&1&0\\1&1&1\end{pmatrix}+ per_{\omega}\begin{pmatrix}1&1&1\\1&0&1\\1&1&1\end{pmatrix})+\omega per_{\omega}\begin{pmatrix}1&1&1\\1&0&0\\1&1&1\end{pmatrix}+  $$
$$\omega^2(\enskip per_{\omega}\begin{pmatrix}1&0&1\\1&1&0\\1&1&1\end{pmatrix}+ per_{\omega}\begin{pmatrix}1&1&0\\1&0&1\\1&1&1\end{pmatrix})+\omega per_{\omega}\begin{pmatrix}1&1&1\\1&1&0\\1&1&1\end{pmatrix}+  $$
$$\omega^2 per_{\omega}\begin{pmatrix}1&1&0\\1&1&1\\1&1&1\end{pmatrix}=(\omega+\omega^2)+(\omega+\omega^2+1)+(2\omega^2+1)+$$
$$(\omega^2+2+\omega)+
(2+2\omega)+(\omega^2+1)+(\omega^2+1+\omega)+(2+\omega)+$$ $$(\omega+2\omega^2+1)+(\omega^2+2+\omega)=13+9\omega+10\omega^2.  $$
Thus in $\mathcal{B}(A)$ we have 13 permutation of class 0 modulo 3; 9 permutations
 of class 1 modulo 3 and 10 permutations of class 2 modulo 3.\newline
$\blacksquare$
\section{On two sequences connected with  $Cycle(A)$ }
In summer of 2010, the author published two sequences A179926 and A180026 in OEIS \cite{10}. $a(n):=A179926(n)$ is defined as the number of permutations of all $\tau(n)$ divisors of $n$ of the form: $d_1=n, \enskip d_2,\enskip d_3,\enskip...,\enskip d_{\tau(n)}$ such that $\frac {d_{i+1}}{d_i}$ is a prime or 1/prime for $i=1,...\tau(n).$ Note that
$a(n)$ is a function of exponents of prime power factorization of $n$ only; moreover, it is invariant with respect to permutations of them. This sequence is equivalently defined as the number of ways, for a given finite multiset $E,$ to get, beginning with $E,$  all
 submultisets of $E,$ if in every step we remove or join one element of $E.$
  Sequence $b(n):=A180026(n)$ differs from A179926 by an additional condition: $\frac {d_{\tau(n)}}{d_1}$ is a prime. In the equivalent formulation it corresponds to
  the condition that in the last step $E$ is obtained from a submultiset
   by joining one element.\newline
\indent Note that, it is easy to prove that, knowing any permissible permutation of divisors, say, $\delta_1=n, \enskip \delta_2,\enskip...,\enskip \delta_{\tau(n)}$ (such that $\frac {\delta_{i+1}}{\delta_i}$ is a prime or 1/prime), we can calculate $b(n)$ using the following construction. Consider square (0,1) matrix $B=\{b_{ij}\}$ of order $\tau(n)$ in which $b_{ij}=1,$ if $\frac {\delta_{i}}{\delta_j}$ is prime or 1/prime, and $b_{ij}=0,$ otherwise. Then $b(n)=Cycle (B).$ In case of A179926, the construction is a little more complicated: $a(n)=Cycle (A),$ where $A$ is obtained from $B$ by the replacing the first its column by the column from
1's.\newpage
\indent Note also that A. Heinz \cite{2} proved that, in particular, $a(\Pi_{i=1}^n p_i),$ where $p_i$ are distinct primes, equals to the number of Hamiltonian paths (or Gray codes) on $n$-cube with a marked starting node (see A003043 in \cite{10}), while $b(\Pi_{i=1}^n p_i)$ equals to the number of directed Hamiltonian cycles on $n$-cube (see A003042 in \cite{10}).
\section{Conclusive remarks}
In spite of the algorithm (\ref{4.8})-(\ref{4.11}) for the evaluation of cycle index of permutations with restricted positions $C(A; t_1, t_2,...,t_n)$ is rather complicated,
it turns out that one can find an explicit representation for it.
 It appears that a unique way for that gives a so-called "method of index of
 arrangements" which was discovered by the author in \cite{6}. This method was
 realized for finding explicit formulas for $C(A; t_1, t_2,...,t_n)$ (with concrete examples) in \cite{9}. Actually, in \cite{9} we gave a wide development of Riordan-Kaplansky 
 theory of rook polynomials \cite{5}.
\section{Acknowledgment}
The author is grateful to Alois Heinz for a computer realization of formula (\ref{4.2}) and the calculations of large values of sequences A179926 and A180026.

\;\;\;\;\;\;\;\;
\end{document}